\documentclass[11pt]{amsart}
\usepackage{amssymb}

\newtheorem{theorem}{Theorem}[section]
\newtheorem{corollary}[theorem]{Corollary}

\newtheorem{proposition}[theorem]{Proposition}

\begin{document}

\title[Hilbertian versus Hilbert W*-modules, $L^2$- and other invariants]{Hilbertian
  versus Hilbert W*-Modules, and applications to $L^2$- and other invariants}
\author[M.~Frank]{Michael Frank}
\address{Universit\"at Leipzig, Mathematisches Institut, Augustusplatz 10,
         D-04109 Leipzig, Fed.~Rep.~Germany}
\email{frank@mathematik.uni-leipzig.de}
\thanks{The research was partially supported by the INTAS grant INTAS 96-1099.}
%\date{April 3, 1997 / November 9, 1999/March 28, 2000}

%%%%%%%%%%%%%%%%%%%%%%%%%%%%%%%%%%%%%%%%%%%%%%%%%%%%%%%%%%%%%%%%%%%%%%%%%%%%%%%
\begin{abstract}
Hilbert(ian) $A$-modules over finite von Neumann algebras with a faithful
normal trace state (from global analysis) and Hilbert W*-modules over $A$
(from operator algebra theory) are compared, and a categorical equivalence
is established. The correspondence between these two structures sheds new
light on basic results in $L^2$-invariant theory providing alternative proofs.
We indicate new invariants for finitely generated projective $B$-modules,
where $B$ is a unital C*-algebra, (usually the full group C*-algebra $C^*(\pi)$
of the fundamental group $\pi=\pi_1(M)$ of a manifold $M$).
\end{abstract}
%%%%%%%%%%%%%%%%%%%%%%%%%%%%%%%%%%%%%%%%%%%%%%%%%%%%%%%%%%%%%%%%%%%%%%%%%%%%%%%

\maketitle

During the last decade W.~L\"uck, A.~Carey, V.~Mathai, and other authors
\cite{Lueck/Rothenberg,Lueck97,Mathai/Carey,Carey/Mathai} used the
analytical concept of Hilbert(ian) $A$-modules over finite von
Neumann algebras $A$ for the study of $L^2$-invariants in global analysis
The technical concept was originally invented by M.~F.~Atiyah \cite{Atiyah} and
I.~M.~Singer \cite{Singer} in 1976, and further developed by J.~Cheeger
and M.~Gromov \cite{Cheeger/Gromov1, Cheeger/Gromov2} in 1985-86. For
the detailed history we refer to the handbook \cite{Lueck98-3}.
The goal has been to obtain
$L^2$-invariants for certain compact closed manifolds, e.g.~$L^2$-(co)homology,
$L^2$-torsion, $L^2$-Betti numbers and Novikov-Shubin invariants for finitely
generated Hilbert $A$-chain complexes. On the other hand Hilbert C*-modules
over arbitrary C*-algebras have been used in operator and operator algebra
theory, in global analysis, in noncommutative geometry and mathematical physics
for about 50 years, \cite{Kaplansky,Lance,NEWO,Frank:97}.
The purpose of the present note is to compare these two categories of
$\mathbb C \,$-/C*-valued inner product modules over finite von Neumann
algebras $A$, where for technical purposes the von Neumann algebras $A$ are
supposed to admit a faithful normal trace state. We establish a categorical
equivalence.
Transferring known results on type ${\rm II}_\infty$ von Neumann algebras and
self-dual Hilbert W*-modules through this categorical equivalence to the
theory of Hilbertian $A$-modules we obtain more evidence on the background of
the theory of $L^2$-invariants from a different viewpoint. We establish new
invariants for finitely generated projective $B$-modules over unital C*-algebras
$B$ and give a perspective for future research.

%%%%%%%%%%%%%%%%%%%%%%%%%%%%%%%%%%%%%%%%%%%%%%%%%%%%%%%%%%%%%%%%%%%%%%%%%%%%%%%
\section{Introduction}

\bigskip
Let $A$ be a finite von Neumann algebra. Since our goal is a comparison of
two categories of Hilbert modules over $A$ which are already described in the
literature we have to suppose the existence of a faithful normal trace state
$tr$ on $A$ to meet the requirements of one of the concepts under
consideration. Basically, this assumption only restricts the center of $A$
with respect to its attainable dimension. The restriction is caused by
the general type and direct integral decomposition theory of von Neumann
algebras, see \cite[V.~Th.~2.4, 2.6]{Takesaki}.
Since finite von Neumann algebras always admit a center-valued trace
functional we can easily extend the comparison of both the categories under
consideration allowing general finite von Neumann algebras. The only loss
would be that we would have to treat Hilbert W*-modules over the center of
finite von Neumann algebras $A$ instead of Hilbert spaces carrying an
additional proper $A$-module structure as described below. However, for
historical reasons we do not invent this more general point of view and
consider only the more restrictive setting supposing the existence of a
trace state on $A$.

\medskip
One way to define a sort of Hilbert(ian) modules over $A$ is the following
often used one which has been introduced in connection with geometric
applications by W.~L\"uck and by A.~Carey, V.~Mathai in 1991/92. Basic
references are the publications by
W.~L\"uck and M.~Rothenberg \cite{Lueck/Rothenberg},
A.~Carey and V.~Mathai \cite{Mathai/Carey,Carey/Mathai,Mathai:98/1,Mathai:98/2},
M.~Rothenberg \cite{Mathai/Rothenberg}, M.~A.~Shubin \cite{Mathai/Shubin:96} and
M.~Farber \cite{Farber,Farber98}, \cite{Carey/Farber/Mathai},
J.~Lott and W.~L\"uck \cite{Lott/Lueck},
D.~Burghelea, L.~Friedlander, T.~Kappeler, P.~McDonald \cite{BFKM,BFK,BFK:99},
W.~L\"uck
\cite{Lueck94-1,Lueck94-2,Lueck94-3,Lueck97,Lueck98-1,Lueck98-2,Lueck98-3},
T.~Schick \cite{Schick,Schick99},
H.~Reich \cite{Reich,LRS:98},
A.~Carey, V.~Mathai and A.~S.~Mishchenko \cite{Carey/Mathai/Mishchenko},
and others. Since the authors differ slightly in their denotations we give a
comprehensive account at this place adopting the denotations introduced by
W.~L\"uck and Th.~Schick in the main.

\smallskip
To introduce Hilbertian modules over finite von Neumann algebras $A$
with trace functional we begin with some {\it standard Hilbert modules over
$A$} -- the completion of the algebraic tensor products $A \odot H$ of $A$ by
Hilbert spaces $H$ with respect to the Hilbert norm induced by the
$\mathbb C$-valued inner product
\[
   \langle a \otimes h, b \otimes g \rangle = tr(ab^*) \cdot \langle h,g
   \rangle_H \, , \,\, {\rm where} \,\, a,b \in A, \, h,g \in H \, ,
\]
The $A$-module action is defined as the (left, w.l.o.g.) multiplication of
$A$ by elements of $A$. So all these sets become Hilbert spaces and (left)
$A$-modules by construction.
Following \cite{Lueck/Rothenberg,Carey/Mathai,Mathai/Carey} they are denoted
by $l^2(A) \otimes H$. A {\it Hilbert module over $A$} is a Hilbert space
$\mathcal M$ together with a continuous (left) $A$-module structure such that
there exists an isometric $A$-linear embedding of $\mathcal M$ into one of the
standard Hilbert modules $l^2(A) \otimes H$ over $A$. A {\it Hilbertian module
over $A$} is a topological vector space $\mathcal M$ with a continuous (left)
$A$-action such that there exists a scalar product $\langle .,. \rangle$ on
$\mathcal M$ generating the topology of $\mathcal M$ and turning $\{
{\mathcal M}, \langle .,. \rangle \}$ together with the $A$-action into a
Hilbert module over $A$ in the previously defined sense. Any two scalar
products $\langle .,. \rangle_1$ and $\langle .,. \rangle_2$ on a Hilbertian
module $\mathcal M$ over $A$ that induce one and the same topology on
$\mathcal M$ and turn $\{ {\mathcal M}, \langle .,. \rangle_1 \}$ and
$\{ {\mathcal M}, \langle .,. \rangle_2 \}$, respectively, into concrete
Hilbert modules over $A$ are linked by a bounded self-adjoint positive
$A$-linear map $T: {\mathcal M} \to \mathcal M$ fulfilling the identity
$\langle .,. \rangle_1 \equiv \langle T(.),T(.) \rangle_2$ on ${\mathcal M}
\times \mathcal M$, \cite[item 1.4]{Carey/Farber/Mathai}.
A Hilbertian module over $A$ is {\it finitely generated} if it is embedable
as a Hilbert module over $A$ into some standard Hilbert module
$l^2(A) \otimes {\mathbb C}^n$ for finite $n \in {\mathbb N}$.

If $\mathcal M$, $\mathcal N$ are two Hilbertian modules over $A$, then a
bounded $A$-linear map $T: {\mathcal M} \to {\mathcal N}$ is called an {\it
$A$-module morphism}. By \cite[I.6.~Th.2]{Dixmier} and \cite[Lemma 9.11]{Lueck97}
the set of all $A$-module morphisms ${\rm End}_A({\mathcal M})$ on a Hilbertian
module $\mathcal M$ has the structure of a von Neumann subalgebra of the
type I von Neumann algebra ${\rm End}_{\mathbb C}({\mathcal M})$ of all bounded
linear operators on the Hilbert space $\mathcal M$. Its commutant is precisely
$A$ itself since the $A$-action on $\mathcal M$ was supposed to be continuous.
Hilbertian modules over a fixed finite von Neumann algebra $A$ with faithful
normal trace state together with the specified $A$-module morphisms form a
category. This category is one of the two we are going to compare.

\bigskip
The second kind of Hilbert modules over C*-algebras is widely used in operator
algebra theory and its applications. There are also topological applications of
this structure, cf.~\cite{Carey/Mathai/Mishchenko,Lueck97,Lueck98-1,Lueck98-2}
for example.
The theory of Hilbert C*-modules goes back to the work of I.~Kaplansky
\cite{Kaplansky}, W.~L.~Paschke \cite{Paschke} and M.~A.~Rieffel \cite{Rieffel}.
Contemporary explanation are published in E.~C.~Lance's lecture notes
\cite{Lance} and in I.~Raeburn's and D.~P.~Williams' book \cite{RaeWil}.
For a detailed bibliography on the subject we refer to
\cite[Appendix]{Frank:97}.

A Hilbert C*-module over a C*-algebra $A$ is a (left, w.~l.~o.~g.) $A$-module
$\mathcal M$ together with a map $\langle .,. \rangle_{{\mathcal M}} :
{\mathcal M} \times {\mathcal M} \to A$ satisfying the following axioms:

\newpage

(i) $\: \langle x,x \rangle_{{\mathcal M}} \geq 0$ for every $x \in \mathcal M$,

(ii) $\, \langle x,x \rangle_{{\mathcal M}} = 0$ if and only if $x=0$,

(iii) $\langle x,y \rangle_{{\mathcal M}} = \langle y,x \rangle_{{\mathcal M}}$
      for every $x,y \in \mathcal M$,

(iv) $\, \langle ax+by,z \rangle_{{\mathcal M}} = a \langle x,z \rangle_{
     {\mathcal M}} + \langle y,z \rangle_{{\mathcal M}}\,$ for every $\,a,b \in
     A$, $x,y,z \in \mathcal M$.

\smallskip \noindent
In this setting the map $\|.\|: {\mathcal M} \to {\mathbb R}^+$ defined by
$\| x \| = \| \langle x,x \rangle_{{\mathcal M}} \|_A^{1/2}$, $(x \in \mathcal
M)$, is a norm on $\mathcal M$. In the sequel we always assume $\mathcal M$ to
be complete with respect to this norm. The pair $\{ {\mathcal M}, \langle .,.
\rangle_{{\mathcal M}} \}$ is said to be a {\it Hilbert $A$-module}, and the
map $\langle .,. \rangle_{{\mathcal M}}$ is referred to as to an {\it
$A$-valued inner product on $\mathcal M$}.

A {\it morphism of Hilbert C*-modules} over a certain C*-algebra $A$ is a
bounded linear map between them intertwining the given actions of $A$ on
the modules. Hilbert C*-modules over a fixed C*-algebra $A$ and bounded
$A$-linear morphisms form another category to be considered.

\smallskip
One of the unpleasant properties of some Hilbert C*-modules is the general
failing of the analog of the Riesz representation theorem for bounded module
maps from the Hilbert C*-module into the C*-algebra of its coefficients,
where the C*-algebra of coefficients is usually assumed to act on itself by
multiplication from the left to become an $A$-module. For example, consider
the finite von Neumann algebra $A = l_\infty$ of all bounded complex sequences
and the Hilbert $A$-module ${\mathcal M}=c_0$ of all sequences converging to
zero equipped with the $A$-valued inner product $\langle x,y \rangle_A = ab^*$,
$(a,b \in A)$.  Every bounded $A$-linear map from $c_0$ to $l_\infty$ is given
by a multiplication of ${\mathcal M}=c_0$ by some element of $l_\infty$ from
the right. However, if this factor belongs to $l_\infty \setminus c_0$ then
the resulting map is not representable as $\langle .,x \rangle_A$ for elements
$x \in c_0$. The next bothering fact is the non-adjointability of some bounded
module operators on Hilbert C*-modules. To present an example, consider
a finite von Neumann algebra $A$ together with one of its maximal left ideals
$I$. Form the Hilbert $A$-module ${\mathcal M} = A \oplus I$ equipped with the
$A$-valued inner product $\langle (a,i),(b,j) \rangle_{{\mathcal M}} =
ab^* + ij^*$, $(a,b \in A, \, i,j \in I)$. Then the operator $T$ defined by
$T((a,i)) = (i,0)$ fails to be adjointable since its formal adjoint operator
does not map $\mathcal M$ back into itself, but partially into a different set.

\smallskip
For Hilbert C*-modules over von Neumann algebras (in short: Hilbert W*-modules)
there is a manner to overcome these obstacles on the way to a 'better' theory
within the category of Hilbert W*-modules. The basic idea is to dilate every
Hilbert W*-module to a larger Hilbert W*-module that has properties analogous
to those known for Hilbert spaces in any respect, and to select a canonically
smallest dilation of that kind. W.~L.~Paschke described the construction of
those envelopes in \cite{Paschke} in 1973. To introduce it we need some further
preparation.

For any C*-algebra $A$ denote the set of all bounded $A$-linear maps from
a given Hilbert $A$-module $\mathcal M$ into $A$ by ${\mathcal M}'$. This set
${\mathcal M}'$ canonically becomes a Banach $A$-module that contains an
isometrically embedded copy of the original Hilbert $A$-module $\mathcal M$
which can be obtained as the image of the map $x \in \mathcal M \to \langle
.,x \rangle \in {\mathcal M}'$ .
In case the C*-algebra $A$ is a von Neumann algebra the $A$-valued inner
product $\langle .,. \rangle_{{\mathcal M}}$ on a Hilbert $A$-module
$\mathcal M \hookrightarrow {\mathcal M}'$ can always be continued to an
$A$-valued inner product on ${\mathcal M}'$ preserving the compatibility with
the standard embedding of $\mathcal M$ into ${\mathcal M}'$,
\cite[Th.~3.2]{Paschke}.
As a result ${\mathcal M}'$ turns out to be a self-dual Hilbert $A$-module
(i.e.~${\mathcal M}' \equiv ({\mathcal M}')'$), and
the Banach algebra of all bounded $A$-linear endomorphisms ${\rm End}_A
({\mathcal M}')$ of ${\mathcal M}'$ becomes a von Neumann algebra,
\cite[Prop.~3.10]{Paschke}. If $A$ is of type I or II or III, respectively,
then ${\rm End}_A({\mathcal M}')$ is always exactly of the same type as $A$,
\cite[Cor.~8.7, Rem.]{Rieffel}.
Moreover, ${\mathcal M}'$ can be constructed from $\mathcal M$ in a simple
topological way: the unit ball of ${\mathcal M}'$ is the completion of
the unit ball of $\mathcal M$ with respect to the topology induced by anyone
of the two sets of semi-norms
\[
   \{ f(\langle .,x \rangle_{{\mathcal M}})^{1/2} \, : \,
   f \in A_*^+, \, \|f\|=1 , \, x \in \mathcal M, \| x \| \leq 1 \}  \, ,
\]
\[
   \{ f(\langle .,. \rangle_{{\mathcal M}})^{1/2} \, : \,
   f \in A_*^+, \, \|f\|=1 \} \, ,
\]
where $A_*^+$ denotes the set of all normal positive functionals on $A$,
\cite{Paschke,BDH}, \cite[Th.~3.2]{Frank:90}. In particular, the biorthogonal
complement of any Hilbert $A$-submodule $\mathcal N \subseteq \mathcal M$ with
respect to the self-dual Hilbert $A$-module ${\mathcal M}'$ can be canonically
identified with the $A$-dual Banach $A$-module ${\mathcal N}'$ since these
topologies are orthogonality-preserving. Among the various fortunate properties
of self-dual Hilbert W*-modules we mention that another $A$-valued inner
product $\langle .,. \rangle_2$ on a self-dual Hilbert $A$-module
$\{ {\mathcal M}, \langle .,. \rangle_1 \}$ induces the same norm topology on
$\mathcal M$ if and only if there exists a bounded $A$-linear positive operator
$T \in {\rm End}_A({\mathcal M})$ such that the identity $\langle .,. \rangle_1
\equiv \langle T(.),T(.) \rangle_2$ holds on ${\mathcal M} \times \mathcal M$.
Furthermore, two self-dual Hilbert $A$-modules are unitarily isomorphic if and
only if they are isomorphic as Banach $A$-modules (\cite[Prop.~2.2]{Frank:90}),
a property which sometimes fails to hold for general Hilbert C*-modules
(cf.~\cite{Frank:99}) and which also causes the slight difference between the
definition of Hilbert C*-modules and that of Hilbert spaces.

One easily checks that every bounded $A$-linear map defined between two
Hilbert $A$-modules over von Neumann algebras $A$ continues to a unique
bounded $A$-linear map between their $A$-dual self-dual Hilbert $A$-modules,
\cite[Cor.~3.7]{Paschke}.
These circumstances make 'self-dualisation' of Hilbert W*-modules a suitable
operation to overcome the arising technical difficulties of the general theory
of Hilbert C*-modules and W*-modules.

Furthermore, self-dual Hilbert W*-modules over von Neumann algebras
$A$ can be alternatively described as w*-closures of formal direct sums
(i.e.~tuples) of *weakly closed left ideals of $A$: every self-dual Hilbert
$A$-module $\{ {\mathcal M}, \langle .,. \rangle_{{\mathcal M}} \}$ is unitarily
isomorphic to a Hilbert $A$-module
\[
{\mathcal N} = \left\{ \{ a_\alpha \}_{\alpha \in I} \, : \,
              a_\alpha \in Ap_\alpha,
              0 \leq p_\alpha = p_\alpha^2 \in A \right\}
             =:  \sum_{\alpha \in I} \oplus \,\, Ap_\alpha \, ,
\]
\[
\langle \{ a_\alpha \}_{\alpha \in I},\{a_\alpha \}_{\alpha \in I} \rangle =
        \sup_{S \in {\mathcal F}}
        \left( \sum_{\alpha \in S} a_\alpha a_\alpha^* \right) \, ,
\]
where $\mathcal F$ denotes the net of all finite subsets of the index set $I$,
\cite[Th.~3.12]{Paschke}.
Obviously, $\mathcal N$ can be characterized isomorphicly as a direct orthogonal
summand of the $A$-dual Hilbert $A$-module ${\mathcal H}'$ of some standard
Hilbert $A$-module ${\mathcal H} = A \otimes H$, where the Hilbert space $H$
possesses a Hilbert basis of cardinality ${\rm card}(I)$. Here the $A$-valued
inner product on $\mathcal H$ is usually defined by $\langle a \otimes h, b
\otimes g \rangle_{\mathcal H} = ab^*\langle h,g \rangle_H$ on elementary
tensors.

%%%%%%%%%%%%%%%%%%%%%%%%%%%%%%%%%%%%%%%%%%%%%%%%%%%%%%%%%%%%%%%%%%%%%%%%%%%%%%%

\section{Hilbertian modules versus Hilbert W*-modules}

\bigskip
Comparing the direct sum decompositions of (standard) Hilbertian and self-dual
Hilbert W*-modules the similarity between these two categories comes to light.
In a first step W.~L\"uck was able to identify the appropriate subcategories
of finitely generated Hilbertian modules over $A$ with the subcategory of
finitely generated projective W*-modules over $A$, \cite[Th.~2.1]{Lueck97},
which are precisely the finitely generated Hilbert $A$-modules by
\cite{Mish,NEWO}.
Our first goal is to establish the hidden categorical identification in
full:

\newpage

\begin{theorem}
  Let $A$ be a finite von Neumann algebra that possesses a normal faithful
  trace state $tr$. The two categories

   (i) $\,$Hilbertian modules over $A$, $A$-module morphisms;

   (ii) self-dual Hilbert W*-modules over $A$, bounded $A$-linear morphisms;

  \noindent
  are equivalent.
  The involution of morphisms is intertwined by the linking functor $\Phi$.
\end{theorem}

\begin{proof}
We construct a functor $\Phi$ from the first category into the second category
and investigate its properties.

Let $\mathcal N$ be a Hilbertian module over $A$. Suppose, $\mathcal N$ is
already $A$-linearly and continuously embedded into a certain standard
Hilbertian module $l^2(A) \otimes H$ for some Hilbert space $H$.
Consider the intersection $\Phi({\mathcal N})$ of $\mathcal N$ with the
norm-dense subset $(A \otimes H)' \subseteq l^2(A) \otimes H$. This intersection
$\Phi({\mathcal N})$ is non-empty and $A$-invariant by construction.
Moreover, for topological reasons it is complete with respect to the norm
$\|.\|=\| \langle .,. \rangle_{(A \otimes H)'} \|_A^{1/2}$, and the unit ball
of $\Phi({\mathcal N})$ is complete with respect to the topology induced by the
semi-norms $\{ f(\langle .,. \rangle_{A \otimes H})^{1/2} : f \in A_*^+ , \,
\|f\|=1 \}$. By \cite[Th.~3.2]{Frank:90} $\Phi({\mathcal N})$ is a self-dual
Hilbert W*-module over $A$, and the completion of it with respect to the norm
$\|.\|=tr(\langle .,. \rangle_{\Phi({\mathcal N})})^{1/2}$ recovers $\mathcal
N$. Note, that $\Phi({\mathcal N})$ is a direct orthogonal summand of the
Hilbert W*-module $(A \otimes H)'$ (\cite[Th.~2.8]{Frank:90}), and that the
$A$-linear projection to it extends to the $A$-linear projection from $l^2(A)
\otimes H$ to $\mathcal N$.

Consider two $A$-linear continuous embeddings of $\mathcal N$ into $l^2(A)
\otimes H_1$ and $l^2(A) \otimes H_2$, respectively. The images of $\mathcal N$
in $l^2(A) \otimes H_1$ and in $l^2(A) \otimes H_2$ are linked by an isometric
$A$-linear operator with carrier projections equal to the projections to
the embedded copies of $\mathcal N$ in $l^2(A) \otimes H_1$ and $l^2(A)
\otimes H_2$, respectively. Repeating our construction we obtain two self-dual
Hilbert W*-modules $\Phi_1({\mathcal N})$ and $\Phi_2({\mathcal N})$
which are $A$-linearly and isometrically isomorphic. By E.~C.~Lance's
theorem \cite[Th.~3.5]{Lance} they are unitarily isomorphic, too.
Consequently, the functor $\Phi$ does not depend on the embedding
of $\mathcal N$ as a Hilbertian module over $A$.

Let $H$ be any Hilbert space. Then $l^2(A) \otimes H$ is the Hilbert norm
closure of the self-dual Hilbert $A$-module $\{ {\mathcal H}=(A \otimes H)',
\langle .,. \rangle_{l_2} \}$ with respect to the Hilbert space norm
$\|.\|=tr(\langle .,. \rangle_{{\mathcal H}})^{1/2}$. For general self-dual
Hilbert $A$-modules $\{ \mathcal H, \langle .,. \rangle_{{\mathcal H}} \}$
consider their canonical decomposition as $\mathcal H = \sum_{\alpha \in I}
\oplus Ap_\alpha$ with $p_\alpha^2=p_\alpha > 0$ of $A$. Obviously, $\mathcal
H$ is an orthogonal direct summand of the self-dual Hilbert $A$-module
$(A \otimes H)'$ for a certain Hilbert space $H$ with $\dim(H)= {\rm card}
(I)$, where the orthogonal complement is the self-dual Hilbert $A$-module
$({\mathcal H})^\bot = \sum_{\alpha \in I} \oplus A(1_A-p_\alpha)$. Hence, the
completion of $\mathcal H$ with respect to the norm $tr(\langle .,.
\rangle_{{\mathcal H}})^{1/2}$ is a Hilbertian module over $A$ by definition.
This shows the functor $\Phi$ to be surjective. At the same time we constructed
the inverse functor $\Phi^{-1}$.

The respective sets of bounded $A$-linear morphisms of both these categories
can be seen to coincide looking at the present construction. Since every
bounded module map on a self-dual Hilbert $A$-module possesses an adjoint
and since every bounded module map on a Hilbertian module $\mathcal N$ over
$A$ preserves the subset $\Phi(\mathcal N)$ invariant, the coincidence of both
the involutions comes to light.
\end{proof}

\begin{corollary}
  The Hilbert space orthogonal complement of a Hilbertian module $\mathcal N$
  over a finite von Neumann algebra $A$ with faithful normal trace state
  inside a standard Hilbertian module $l^2(A) \otimes H$ is $A$-invariant
  and, hence, a Hilbertian module, too.
\end{corollary}

\begin{corollary} {\rm (W.~L\"uck's theorem, \cite[Th.~2.1]{Lueck97}),
    \cite[Th.~1.8]{Lueck98-1}} \newline
  The functor $\Phi$ identifies the subcategory of finitely generated
  Hilbertian modules over $A$ with the subcategory of finitely generated
  projective $A$-modules (i.e.~finitely generated Hilbert $A$-modules).
\end{corollary}

For a proof we have only to recall that $\Phi(l^2(A) \otimes {\mathbb C}^N) =
(A \otimes {\mathbb C}^N)' = A \otimes {\mathbb C}^N$ for every finite
$N \in {\mathbb N}$.

\begin{corollary}
  Let $\mathcal N$ be a Hilbertian module over a certain finite von Neumann
  algebra $A$ that admits a normal trace functional $tr$. Suppose, $\mathcal N$
  can be identified with a direct summand of the standard Hilbertian module
  $l^2(A) \otimes H$ for a certain Hilbert space $H$.
  Then the von Neumann algebra ${\rm End}_A(l^2(A) \otimes H)$ of all bounded
  $A$-linear module maps on $l^2(A) \otimes H$ is $*$-isomorphic to the
  W*-tensor product $A \,\overline{\otimes}\, {\rm End}_{\mathbb C}(H)$ and, hence,
  a type ${\rm II}_\infty$ von Neumann algebra. The von Neumann algebra
  ${\rm End}_A(\mathcal N)$ can be identified with a full corner of the von
  Neumann algebra ${\rm End}_A(l^2(A) \otimes H)$, i.e. $P \cdot {\rm End}_A
  (l^2(A) \otimes H) \cdot P \equiv {\rm End}_A(\mathcal N)$ for the orthogonal
  projection $P: l^2(A) \otimes H \to \mathcal N$.

  Fixing the trace functionals $tr$ and $Tr_{B(H)}$ on $A$ and $B(H)$,
  respectively, the canonical semifinite trace functional $Tr$ on ${\rm End}_A
  (l^2(A) \otimes H)$ defined by $Tr(a \otimes T) = tr(a) \cdot Tr_{B(H)} (T)$
  on elementary tensors assigns either a real number $Tr(P) \in [0,+\infty)$
  or the symbol $+\infty$ to $\mathcal N$. The assigned value does not depend on
  the choice of the embedding of $\mathcal N$, and the assigned value is a finite
  number if and only if $\mathcal N$ is finitely generated.
\end{corollary}

These elementary conclusions can be derived from standard von Neumann algebra
theory identifying the Hilbertian module $\mathcal N$ with its image $\Phi(
\mathcal N)$ in the category of Hilbert W*-modules over $A$, cf.~\cite{Paschke}.
If we have two embeddings of $\mathcal N$ into standard Hilbertian modules
$l^2(A) \otimes H_1$ and $l^2(A) \otimes H_2$, respectively, we may always
assume that both $H_1$ and $H_2$ are isomorphic Hilbert spaces enlarging the
smaller one appropriately. Since there exists a partial isometry $U$ between
the orthogonal projections $P_1$ on $H_1$ and $P_2$ on $H_2$ we obtain
$Tr(P_2) = Tr(UP_2U^*) = Tr(P_1)$, and so the value is independent of the
concrete representation of $\mathcal N$.

%%%%%%%%%%%%%%%%%%%%%%%%%%%%%%%%%%%%%%%%%%%%%%%%%%%%%%%%%%%%%%%%%%%%%%%%%%%%%%%
\section{$L^2$- and other invariants of C*-modules -- revisited}

\bigskip
The purpose of the present section is to investigate the freedom of choice
for (center-valued) semifinite trace functionals on the set of all bounded
module operators on some Hilbertian modules, and to translate the obtained
results into the language of $L^2$-invariants. We rely on the categorical
equivalence obtained in the previous section. So the results are much easier
to obtain for self-dual Hilbert W*-modules because of their well-known
composition structure.
As a result we obtain that the consideration of finitely generated Hilbert
$B$-modules and their invariants can replace the investigation of the
corresponding finitely generated Hilbert $B^{**}$-modules which arise as
categorically equivalent objects of the w*-completion of the former with
respect to the respective $\mathbb C$-valued inner product $tr(\langle .,.
\rangle)$.

In applications of Hilbertian modules ${\mathcal M} \subseteq l^2(A) \otimes H$
over finite von Neumann algebras $A$ with faithful normal trace state $tr$
in the theory of $L^2$-invariants the properties of the canonical projection
$P_{\mathcal M}^H : l^2(A) \otimes H \to \mathcal M$ and the stability of some
of these properties with respect to changes of the representations of
$\mathcal M$ are of major interest, \cite{Lueck97,Carey/Farber/Mathai}. One of
these invariant properties of the family $\{ P_{\mathcal M}^H : H  -  {\rm
Hilbert} \: {\rm space} \}$ is the (non-)existence of a finite ${\rm Z}
(A)$-valued trace value for it, where $Z(A)$ denotes the center of $A$ which
can be canonically identified with the center of ${\rm End}_A(l^2(A) \otimes
H)$. Another property is the existence of a Fuglede-Kadison determinant for
positive invertible module operators on a Hilbertian $A$-module.
If a finite trace value exists then generalized Betti-numbers $b_p^u(C)$ of
chain complexes $C$ of Hilbertian modules over $A$ and Novikov-Shubin invariants
$\alpha_p(C)$ can be defined, \cite{Lueck97}.
However, the value of the trace in $Z(A)$ is not uniquely determined,
cf.~\cite[V.~Th.~2.34]{Takesaki} for a description of the variety of faithful
semifinite normal extended center-valued traces on type ${\rm II}_\infty$ von
Neumann algebras. To fix a standard ${\rm Z}(A)$-valued trace $\tau_o$ on a
C*-subalgebra of ${\rm End}_A(l^2(A) \otimes H)$ we use a standard
$*$-isomorphism
\[
  {\rm End}_A(l^2(A) \otimes H)
        \cong A \, \overline{\otimes} \, {\rm End}_{{\mathbb C}}(H)
\]
(where $\overline{\otimes}$ denotes the W*-tensor product of W*-algebras),
and we set
\[
     \tau_o(a \otimes S) = \tau(a) \cdot Tr_{B(H)}(S)
                     = \tau \otimes Tr_{B(H)}) (a \otimes S)
\]
for $a \in A$ and trace class operators $S$ on $H$. Here $\tau$ denotes the
unique faithful normal center-valued trace on $A$
(cf.~\cite[V.~Th.~2.6]{Takesaki}), and $Tr_{B(H)}$ is a normal trace on $B(H)$
normalized by the requirement that it takes the value one on projections to
one-dimensional subspaces of $H$, \cite{Albematt}.
The following theorem gives a precise description of the situation.
It was already partially stated by W.~L\"uck in \cite[Cor.~3.2]{Lueck97}. In
difference to his result we do neither require the existence of a finite
faithful trace functional on the finite von Neumann algebra $A$ nor have any
preference for any specific center-valued faithful trace on the C*-algebra of
all bounded module operators. The standard trace $\tau_o$ is only a tool to
prove the assertions.

\begin{theorem}
Let $A$ be a finite von Neumann algebra and $\mathcal M$ be a self-dual
Hilbert $A$-module.
Then $\mathcal M$ is finitely generated if and only if for some/every
isometric embedding of $\mathcal M$ into a standard self-dual Hilbert
$A$-module $(A \otimes H)'$ the canonical projection $P_{\mathcal M}^H$ from
$(A \otimes H)'$ to $\mathcal M$ possesses a finite center-valued trace value
$\tau_1(P_{\mathcal M}^H)$, where $\tau_1$ denotes an arbitrary faithful
semifinite normal extended center-valued trace given on the type
${\rm II}_\infty$ von Neumann algebra ${\rm End}_A((A \otimes H)')$.
In other words, the value $\tau_1(P)$ of a center-valued trace $\tau_1$
applied to an orthogonal projection $P \in {\rm End}_A((A \otimes H)')$ is
finite if and only if $P((A \otimes H)')$ is a finitely generated Hilbert
$A$-module, despite the concrete value depends on the choice of $\tau_1$.
\end{theorem}

\begin{proof}
Every isometric copy of $\mathcal M$ as an $A$-submodule of another Hilbert
$A$-module $\mathcal N$ is an orthogonal summand of $\mathcal N$ since
$\mathcal M$ is self-dual by assumption, cf.~\cite{Frank:90}.
In case ${\mathcal M} = (A \otimes H)'$ for some Hilbert space $H$ we have
the isometric algebraic embedding $A \odot {\rm K}_{\mathbb C}(H) \subseteq
{\rm End}_A((A \otimes H)')$, where $\odot$ denotes the algebraic tensor
product. Let us fix the standard center-valued trace $\tau_o$ on $A \odot
{\rm K}_{{\mathbb C}}(H)$. The bounded module operators on $\mathcal M =
(A \otimes H)'$ that admit a finite center-valued trace are contained in the
C*-algebra of all 'compact' operators on $\mathcal M$ which is defined as
the norm-closure of the linear hull of the operators $\{
\theta_{x,y} : \theta_{x,y}(z) = \langle z,x \rangle y \:\, {\rm for} \:\,
x,y,z \in \mathcal M \}$, since the trace class operators $S \in B(H)$ are
all compact operators.
Therefore, projections onto a Hilbert C*-module $\mathcal M = (A \otimes H)'$
admit a finite value with respect to the fixed standard semifinite center-valued
trace $\tau_o$ if and only if they are 'compact' operators, i.e.~if and only
if their image is a finitely generated Hilbert C*-module,
\cite[Th.~15.4.2, Remark 15.4.3]{NEWO}.

By \cite[Th.~2.34]{Takesaki} the type ${\rm II}_\infty$ von Neumann algebra
${\rm End}_A((A \otimes H)')$ may admit more then just one faithful semifinite
normal extended center-valued trace, and an easy classification is available.
However, if a projection $P \in {\rm End}_A((A \otimes H)')$ admits a finite
center-valued trace value with respect to a certain such trace on ${\rm End}_A
((A \otimes H)')$, then the von Neumann algebra $P \cdot {\rm End}_A
((A \otimes H)') \cdot P$ is of type ${\rm II}_1$, i.e.~also the restriction
of the standard faithful center-valued trace $\tau_o$ to it would give only
finite values in $Z(A)$. The converse obviously also holds. This implies our
statement.
\end{proof}

The theorem has some consequences for the point of view taken for the
investigations in the field of $L^2$- and other invariants. The usual approach
is to consider the (full) group C*-algebra $C^*(\pi)$ of the fundamental group
$\pi=\pi_1(M)$ of a certain compact manifold $M$ or, more specifically, of a
finite connected CW-complex $M$. The C*-algebra $C^*(\pi)$ admits a canonical
finite trace functional, so its bidual von Neumann algebra $\mathcal N(\pi)$
has a finite normal trace functional and, hence, is the block-diagonal sum
of $\sigma$-finite type ${\rm II}_1$ and finite type I components. The theorem
above suggests that it might be more perspective to consider the category of
finitely generated projective $C^*(\pi)$-modules instead of the category of
finitely generated projective $\mathcal N(\pi)$-modules. This thought is
supported by a number of results of A.~Carey, V.~Mathai and A.~S.~Mishchenko
\cite{Carey/Mathai/Mishchenko} and of D.~Burghelea, L.~Friedlander and
T.~Kappeler \cite{BFK:99}. They consider the analytic torsion of cone complexes
that arise from finite-dimensional non-simply connected Riemannian manifolds
$M$ and their de Rham complexes $\Omega^*(M;C^*(\pi))$. Their investigations
rely merely on the properties of $C^*(\pi)$ and of finitely generated projective
$C^*(\pi)$-modules, and they avoid any weak* completions of appearing structural
elements of the basic constructions. So one might get the idea that
Betti-numbers and Novikov-Shubin invariants could be already derived in case
$C^*(\pi)$ admits a finite center-valued trace. However, unitarily non-isomorphic
finitely generated projective $C^*(\pi)$-modules can have equal dimension values
in case $C^*(\pi)$ is infinite-dimensional. Furthermore, the monoid of all
finitely generated projective $C^*(\pi)$-modules can fail to have the cancellation
property. Also, the extension of the center-valued trace on $C^*(\pi)$ to the
standard semifinite center-valued trace on ${\rm End}_A(l^2(A) \otimes l^2)$
may assign a finite value to some orthogonal projection $P$ even if the image
$P(l^2(A) \otimes l^2)$ is not finitely generated (cf.~G.~G.~Kasparov's theorem
\cite[Th.~6.2]{Lance} applied to countably generated ideals of $C^*(\pi)$).
So the entire approach does not work for the situation of general finitely
generated projective $C^*(\pi)$-modules without completing them with respect to
the standard weak topology.

\smallskip
However, there is another way to classify finitely generated $B$-modules over
arbitrary unital C*-algebras $B$. It is suggested by the modular frame theory
of countably generated Hilbert $B$-modules that has been recently worked out
by D.~R.~Larson and the author in \cite{FL,FL:99}. Let $B=C^*(\pi)$ and consider
a finitely generated projective $A$-module $\mathcal M$. The module $\mathcal M$
can be equipped with an $B$-valued inner product $\langle .,. \rangle$ and
becomes a finitely generated Hilbert $B$-module that way. Conversely, every
finitely generated Hilbert $B$-module is projective as an $B$-module,
\cite[Cor.~15.4.8]{NEWO}. The choice of the $B$-valued inner product on
$\mathcal M$ is unique up to a positive invertible bounded module operator $T$
on $\mathcal M$ linking any other inner product structure to a fixed one.

By \cite{FL,FL:99} every finitely generated Hilbert $B$-module $\{ \mathcal M,
\langle .,. \rangle \}$ admits at least one finite normalized tight (modular)
frame, i.e.~a $k$-tuple $\{ x_1, ..., x_k \}$ of elements of $\mathcal M$ such
that the equality $\langle x,x \rangle = \sum_{i=1}^k \langle x,x_i \rangle
\langle x_i,x \rangle$ holds for any $x \in \mathcal M$. Furthermore, a
reconstruction formula $x=\sum_{i=1}^k \langle x,x_i \rangle x_i$ is valid for
any $x \in \mathcal M$, so the knowledge of the values $\{ \langle x_i,x_j \}
: i,j=1,...,k \}$ turns out to be sufficient to describe the $B$-module
$\mathcal M$ up to uniqueness. Note that the elements $\{ x_1,...,x_k \}$ need
not to be ($B$-)linearly independent, in general.

\begin{theorem}
  Let $B$ be a unital C*-algebra and let $\{ \mathcal M, \langle .,.
  \rangle_{\mathcal M} \}$ and $\{ \mathcal N, \langle .,. \rangle_{\mathcal N}
  \}$ be two finitely generated Hilbert $B$-modules. Then the following
  conditions are equivalent:

  \newcounter{marke}
   \begin{list}{(\roman{marke})}{\usecounter{marke}}
    \item  $\mathcal M$ and $\mathcal N$ are algebraically isomorphic as
       projective $B$-modules.
    \item  $\{ \mathcal M, \|.\|_{\mathcal M} \}$ and $\{ \mathcal N,
       \|.\|_{\mathcal N} \}$ are isometrically isomorphic as Banach
       $B$-modules.
    \item $\{ \mathcal M, \langle .,. \rangle_{\mathcal M}$ and $\{ \mathcal
       N, \langle .,. \rangle_{\mathcal N} \}$ are unitarily isomorphic as
       Hilbert $B$-modules.
    \item There are finite normalized tight (modular) frames $\{ x_1,...,x_k
       \}$ and $\{ y_1,...,y_l \}$ of $\mathcal M$ and $\mathcal N$,
       respectively, such that $k=l$, $x_i \not= 0$ and $y_i \not= 0$ for any
       $i=1,...,k$, and $\langle x_i,x_j \rangle_{\mathcal M} = \langle y_i,y_j
       \rangle_{\mathcal N}$ for any $i,j=1,...,k$.
   \end{list}
\end{theorem}

\begin{proof} The equivalence of the conditions (i), (ii) and (iii) has been
shown for countably generated Hilbert $B$-modules in \cite[Th.~4.1]{Frank:99}.
The implication (iii)$\to$(iv) can be shown to hold setting $y_i = U(x_i)$
for the existing unitary operator $U: \mathcal M \to \mathcal N$ and for
$i=1,...,k$. The demonstration of the inverse implication requires slightly
more work. For the given normalized tight (modular) frames $\{ x_1,...,x_k \}$
and $\{ y_1,...,y_l \}$ of $\mathcal M$ and $\mathcal N$, respectively, we
define a $B$-linear operator $V$ by the rule $V(x_i)=y_i$, $i=1,...,k$.
For this operator $V$ we have
  \begin{eqnarray*}
     \langle V(x),y_j \rangle_{\mathcal N} & = &
        \sum_{i=1}^k \langle V(x),y_i \rangle_{\mathcal N} \langle y_i,y_j
        \rangle_{\mathcal N} \\
     &=& \sum_{i=1}^k \left\langle \sum_{m=1}^k \langle x,x_m \rangle_{\mathcal M}
        V(x_m),y_i \right\rangle_{\mathcal N} \langle x_i,x_j \rangle_{\mathcal M}
  \end{eqnarray*}
  \begin{eqnarray*}
     \quad\quad
     &=& \sum_{i=1}^k \sum_{m=1}^k \langle x,x_m \rangle_{\mathcal M} \langle x_m,
        x_i \rangle_{\mathcal M} \langle x_i,x_j \rangle_{\mathcal M} \\
     &=& \sum_{m=1}^k  \langle x,x_m \rangle_{\mathcal M} \left\langle x_m,
        \sum_{i=1}^k \langle x_j,x_i \rangle_{\mathcal M} x_i
        \right\rangle_{\mathcal M} \\
     &=& \left\langle x , \sum_{m=1}^k \langle x_j,x_m \rangle_{\mathcal M} x_m
        \right\rangle_{\mathcal M} \\
     &=& \langle x,x_j \rangle_{\mathcal M}
  \end{eqnarray*}
for every $x \in \mathcal M$. Consequently,
\[
  \langle x,x \rangle_{\mathcal M} =
  \sum_{i=1}^k \langle x,x_i \rangle_{\mathcal M} \langle x_i,x \rangle_{\mathcal M}=
  \sum_{i=1}^k \langle V(x),y_i \rangle_{\mathcal N} \langle y_i,V(x) \rangle_{\mathcal N}=
  \langle V(x),V(x) \rangle_{\mathcal N}
\]
for any $x \in \mathcal M$, and the operator $V$ is unitary.
\end{proof}

\begin{corollary}
Every finitely generated projective $B$-module $\mathcal M$ over a unital
C*-algebra $B$ can be reconstructed up to isomorphism from the following data:
   \begin{list}{(\roman{marke})}{\usecounter{marke}}
    \item   A finite set of algebraic non-zero modular generators $\{
        x_1,...,x_k \}$ of $\mathcal M$.
    \item   A symmetric $k \times k$ matrix $( a_{ij} )$ of elements from $B$,
        where $a_{ij}$ is supposed to be equal to $\langle x_i,x_j \rangle_0$
        for $i,j = 1,...,k$ and for the (existing and unique) $B$-valued inner
        product $\langle .,. \rangle_0$ on $\mathcal M$ that turns the set of
        algebraic modular generators $\{ x_1,...,x_k\}$ into a normalized
        tight modular frame of the Hilbert $B$-module $\{ \mathcal M,
        \langle .,. \rangle_0 \}$.
   \end{list}
The number of elements in sets of algebraic modular generators of $\mathcal M$
has a minimum, and it suffices to consider sets of generators of minimal
length. Then the modular invariants can be easier compared permuting the
elements of the generating sets if necessary.
\end{corollary}

\begin{proof}
The set of algebraic generators $\{ x_1,...,x_k \}$ of $\mathcal M$ is a frame
with respect to any $B$-valued inner product on $\mathcal M$ which turns
$\mathcal M$ into a Hilbert $B$-module. That is the inequality
\[
  C \cdot  \langle x,x \rangle \leq \sum_{i=1}^k \langle x,x_i \rangle
  \langle x_i,x \rangle \leq D \cdot \langle x,x \rangle
\]
is satisfied for two finite positive real constants $C,D$ and any $x \in
\mathcal M$, see \cite[Th.~5.9]{FL}. What is more, for any frame of $\mathcal
M$ there exists another $B$-valued inner product $\langle .,. \rangle_0$ on
$\mathcal M$ with respect to which it becomes normalized tight, that is
$\langle x,x \rangle_0 = \sum_{i=1}^k \langle x,x_i \rangle_0 \langle x_i,x
\rangle_0$ holds for any $x \in \mathcal M$. This inner product is unique,
cf.~\cite[Cor.~4.3, Th.~6.1]{FL}, \cite[Th.~4.4]{FL:99}.
So assertion (iv) of the previous theorem gives the complete statement.
\end{proof}

\smallskip
We can say more in case the finitely generated Hilbert $B$-module contains
a modular Riesz basis, i.e.~a finite set of modular generators $\{ x_1,...,x_k
\}$ such that the equality $0= b_1x_1 + ... + b_kx_k$ holds for certain
coefficients $\{ b_1,...,b_k \} \subset B$ if and only if $b_ix_i=0$ for any
$i=1,...,k$. Obviously, a modular Riesz basis is minimal as a set of modular
generators, i.e.~we cannot drop any of its elements preserving the generating
property. However, there can exist totally different Riesz bases for the same
module that consist of less elements, cf.~\cite[Ex.~1.1]{FL}. Note
that the coefficients $\{ b_1,...,b_k \}$ can be non-trivial even if $b_ix_i=0$
for any index $i$ since every non-trivial C*-algebra $B$ contains zero-divisors.
Not every Hilbert C*-module containing a normalized tight modular frame does
possess a modular Riesz basis. An example can be found in \cite[Ex.~2.4]{FL:99}.

In case of finitely generated projective W*-modules (and therefore, in the case
of Hilbertian modules over finite W*-algebras) we are in the pleasant situation
that they always contain a modular Riesz basis by \cite[Th.~3.12]{Paschke}.
Moreover, by spectral decomposition every element $x$ of a Hilbert W*-module
$\mathcal M$ has a carrier projection of $\langle x,x \rangle$ contained in the
W*-algebra of coefficients $B$. So we can ascertain the following fact:

\begin{proposition}
  Let $\mathcal M$ be a finitely generated projective $B$-module over a
  W*-algebra $B$ that possesses two finite modular Riesz bases $\{ x_1,...,x_k
  \}$ and $\{y_1,...,y_l \}$.
  Then there exists an $l \times k$ matrix $C = (c_{ij})$, $i=1,...,l$,
  $,j=1,...,k$, with entries from $B$ such that $y_i = \sum_{j=1}^k c_{ij} x_j$
  for any $i = 1,...,l$, and analogously, there exists a $k \times l$ matrix
  $D = (d_{ji})$ with entries from $B$ such that $x_j = \sum_{i=1}^l d_{ji} y_i$
  for any $j=1,...,k$.

  Suppose the left carrier projections of $c_{ij}$ and $d_{ji}$ equal the
  carrier projections of $\langle y_i,y_i \rangle$ and $\langle x_j,x_j \rangle$,
  respectively, and the right carrier projection of $c_{ij}$ and $d_{ji}$ equal
  the carrier projections of $\langle x_j,x_j \rangle$ and $\langle y_i,y_i
  \rangle$, respectively. Then the matrices $C$ and $D$ are Moore-Penrose
  invertible in $M_{kl}(B)$ and $M_{lk}(B)$, respectively. The matrix $C$ is
  the Moore-Penrose inverse of $D$, and vice versa.
\end{proposition}

\begin{proof}
Since both the modular Riesz bases are sets of modular generators of $\mathcal M$
we obtain two $B$-valued (rectangular, w.l.o.g.) matrices $C = (c_{ij})$ and
$D = (d_{ji})$ with $i=1,...,l$ and $j=1,...,k$ such that
\[
  y_i = \sum_{m=1}^k c_{im}x_m \quad , \quad  x_j = \sum_{n=1}^l d_{jn} y_n
   \, .
\]
Combining these two sets of equalities in both the possible ways we obtain
\[
  y_i = \sum_{n=1}^l \left( \sum_{m=1}^k c_{im} d_{mn} \right) y_n  \quad , \quad
  x_j = \sum_{m=1}^k \left( \sum_{n=1}^l d_{jn} c_{nm} \right) x_m
\]
for $i=1,...,l$, $j=1,...,k$. Now, since we deal with sets of coefficients
$\{ c_{ij} \}$ and $\{ d_{ji} \}$ that are supposed to admit special carrier
projections,  the coefficients in front of the elements $\{ y_n \}$ and $\{
x_m \}$ at the right side can only take very specific values:
\[
   \sum_{m=1}^k c_{im} d_{mn} = \delta_{in} \cdot q_n \quad , \quad
   \sum_{n=1}^l d_{jn} c_{nm} = \delta_{jm} \cdot p_m  \, ,
\]
where $\delta_{ij}$ is the Kronecker symbol, $p_m \in B$ is the carrier
projection of $\langle x_m,x_m \rangle$ and $q_n \in B$ is the carrier
projection of $\langle y_n,y_n \rangle$. So $C \cdot D$ and $D \cdot C$ are
positive idempotent diagonal matrices with entries from $B$. The Moore-Penrose
relations $C \cdot D \cdot C = C$, $D \cdot C \cdot D = D$, $(C \cdot D)^* =
C \cdot D$ and $(D \cdot C)^* = D \cdot C$ turn out to be fulfilled.
\end{proof}

\medskip
Summing up, we can replace the single center-valued trace value that
characterizes a finitely generated projective $A$-module over a finite von
Neumann algebra $A$ up to isomorphism by another set of data of a finitely
generated projective $B$-module over a unital C*-algebra $B$. It consists
of a set of non-zero modular generators $\{ x_1, ..., x_k \}$ together with
a $k \times k$ matrix with entries $\langle x_i,x_j \rangle_0$, where $\langle
.,. \rangle_0$ is supposed to be the $B$-valued inner product on the module
with respect to which the generating set becomes a normalized tight frame.
In the von Neumann case the generators can be described as $k$-tuples
$x_i = (0,...,0,p_i,0,...,0)$ for orthogonal projections $p_i \in A$,
$i=1,...,k$, and the $k \times k$ matrix can be chosen to be a diagonal one
with the elements $p_i \in A$ on the diagonal. This follows from the general
structure of self-dual Hilbert W*-modules as described by W.~L.~Paschke at
\cite[Th.~3.12]{Paschke}. For general finitely generated Hilbert $B$-modules
a diagonal structure of the $k \times k$ matrix may not exists for any $k \in
\mathbb N$ since cancellation may not hold in the monoid of all finitely
generated projective $B$-modules. We would like to formulate the problem
whether Moore-Penrose type transfer matrices between modular Riesz
bases appear for more general C*-algebras of coefficients then W*-algebras
or monotone complete C*-algebras, or not.

\smallskip
To derive appropriate invariants for finitely generated Hilbert $B$-chain
complexes in case $B$ is the full group C*-algebra of the fundamental group
of a certain compact manifold $M$ has to await another time since the
considerations would go beyond the set limits of the present paper.

\medskip \noindent
{\bf Acknowledgement:} The author is grateful to V.~M.~Manuilov,
A.~S.~Mishchenko, G.~K.~Pedersen and E.~V.~Troitsky for the fruitful discussions,
the exchange of ideas and the continuous support during the recent years of
collaboration.

%%%%%%%%%%%%%%%%%%%%%%%%%%%%%%%%%%%%%%%%%%%%%%%%%%%%%%%%%%%%%%%%%%%%%%%%%%%%%

\end{document}